%
%
%
%

\input amstex
\documentstyle{amsppt}
\def\bold{\Bbb}
\def\phi{\varphi}
\NoBlackBoxes
\topmatter

\title
The Kobayashi metric for non-convex complex ellipsoids
\endtitle

\author
Peter Pflug and W\l odzimierz Zwonek
\endauthor

\address
Universit\"at Osnabr\"uck Standort Vechta,
Fachbereich Naturwissenschaften/Ma\-the\-matik, Postfach 15 53, D-49364 Vechta
Germany
\endaddress

\email
pflugvec at dosuni1.rz.uni-osnabrueck.de
\endemail
\address
Uniwersytet Jagiello\'nski, Instytut Matematyki, Reymonta 4, 30-059 Krak\'ow,
Poland
\endaddress
\address
current address: Universit\"at Osnabr\"uck Standort Vechta,
Fachbereich Naturwissenschaften/Mathematik
Postfach 15 53, D-49364 Vechta Germany
\endaddress
\email
zwonek at im.uj.edu.pl
\endemail

\abstract
In the paper we give some necessary conditions for a mapping to be
a
$\kappa$-geodesic in non-convex complex ellipsoids. Using these
results we
calculate explicitly the Kobayashi metric in the ellipsoids
$\{|z_1|^2+
|z_2|^{2m}<1\}\subset\bold C^2$, where $m<\frac12$.
\endabstract
\endtopmatter

\subheading{1.Introduction} By $E$ we denote the unit disk in
$\bold C$. Let
$D$ be a domain in $\bold C^n$. For $(z,X)\in D\times \bold C^n$
we
define
$$
\kappa_D(z;X)=\inf \{\gamma_E(\lambda;\alpha):\exists \phi
:E\longrightarrow D \text{ holomorphic, }\phi (\lambda)=z,\alpha
\phi^{\prime}(\lambda)=X\},\tag{1}
$$
where $\gamma_E(\lambda;\alpha)=\frac{|\alpha|}{1-|\lambda|^2}$.
The function
$\kappa_D$ is called \it the Kobayashi pseudometric\rm.
\par We say that a holomorphic mapping $\phi:E\longrightarrow D$
is \it a
$\kappa_D$-geodesic for $(z,X)$ \rm if $\phi$ is a function
achieving the
minimum in \thetag{1}. Without loss of generality we may assume
that
$\phi(0)=z$. Similarly one has a notion of
\it $k$-geodesic \rm for a pair of points, where $k$ denotes the
Kobayashi
pseudodistance (cf. \cite{JP}). We know that if $D$ is a taut
domain then for any
$(z,X)\in D\times\bold C^n$ there is a $\kappa_D$-geodesic for
$(z,X)$ (see
\cite{JP}).
\par If $D$ is a convex domain, then any $\kappa_D$-geodesic $\phi$
for $(z,X)$
with $X\neq 0$ is a $\kappa_D$-geodesic for any
$(\phi(\lambda),\phi^{\prime}(\lambda))$,
$\lambda\in E$; moreover, the Kobayashi pseudometric and the
Carath\'eodory
pseudometric on $D$ coincide (see \cite{L}). In this case any
$\kappa_D$-geodesic for some $(z,X)$ with $X\neq 0$ is called \it
a complex
geodesic \rm (see \cite{Ves}).
\par Let us define
$$
\Cal E(p):=\{ |z_1|^{2p_1}+\ldots + |z_n|^{2p_n}<1\}
\subset \bold C^n,
$$
where $n>1$ and $p=(p_1,\ldots ,p_n)$ with $p_j>0$. We call
$\Cal E(p)$ \it a complex ellipsoid\rm. It is well known that
complex
ellipsoids are
taut domains and they are convex iff $p_j\geq \frac12$ for
$j=1,\ldots,n$.
Moreover, $\partial\Cal E(p)$ is $C^{\omega}$ and strongly
pseudoconvex
at all points $z\in\left(\partial \Cal E(p)\right)\cap(\bold
C_*)^n$.
\par Let us define the following mappings
$\phi=(\phi_1,\ldots,\phi_n)
:E\longrightarrow \bold C^n$, which will be crucial for our
considerations:
$$
\phi_j(\lambda)=\left (\frac{\lambda -\alpha_j}{1-\bar
\alpha_j\lambda}
\right )^{r_j}
\left
(a_j\frac{1-\bar\alpha_j \lambda}{1-\bar \alpha_0\lambda}\right )
^\frac{1}{p_j}\quad (j=1,\ldots,n)
\tag{2}
$$
fulfilling the following conditions:
$$
\gather
a_j\in \bold C_*,\alpha_j\in \bar E\text{ for $j=1,\ldots,n$ and }
\alpha_0\in E,\tag{3}\\
r_j\in\{0,1\}\text{ for $j=1,\ldots,n$ and if $r_j=1$, then
$\alpha_j\in E$,}
\tag{4}\\
\alpha_0=\sum_{j=1}^{n} |a_j|^{2}\alpha_j,\tag{5}\\
1+|\alpha_0|^2=\sum_{j=1}^{n}|a_j|^{2}(1+|\alpha_j|^2).\tag{6}
\endgather
$$
One can easily check that if $\phi$ is not a constant mapping,
then
$\phi(E)\subset \Cal E(p)$.
\par The mappings given by the formulas above  turn out to be the
complex geodesics
in convex complex ellipsoids (certainly under the condition that
they are
not constant). To simplify our notation we will speak in the
following of
$\kappa$-geodesics instead of $\kappa_{\Cal E(p)}$-geodesics
in the case of complex ellipsoids.
\proclaim\nofrills{Theorem 1}\;{\rm (see \cite{JPZ}).} A
non-constant, bounded
mapping $\phi=(\phi_1,\ldots,\phi_n):
E\longrightarrow \bold C^n$, where $\phi_j\not\equiv 0$ for
$j=1,\ldots,n$,
is a complex geodesic in the convex ellipsoid $\Cal E(p)$ iff
$\phi$
is of the form as in \thetag{2} with \thetag{3}--\thetag{6}.
\endproclaim
\par Observe that the assumption $\phi_j\not\equiv 0$ does not
restrict
the generality of this result, since mappings with $0$-components
can be
thought as mappings into a lower dimensional situation.
\par The mappings of \thetag{2} are well-defined not only for
convex
ellipsoids but also for the non-convex ones. Therefore the
natural question arises what these formulas represent in the
general, not
necessarily convex case. Below we shall deal with this problem.
The results we get suggest that all $\kappa$-geodesics are
necessarily of the
form \thetag{2} with \thetag{3}--\thetag{6}.

\proclaim{Proposition 2} Let $\phi:E\longrightarrow \Cal E(p)$ be
a
$\kappa$-geodesic for $(\phi(0),\phi^{\prime}(0))$ with
$\phi^{\prime}(0)\neq
0$, where $\phi_j\not\equiv 0$ for $j=1,\ldots,n$.
Then
$$
\phi_j(\lambda)=B_j(\lambda)\left(a_j\frac{1-\bar\alpha_j\lambda}{1-\bar
\alpha_0\lambda}\right)^{\frac{1}{p_j}},
$$
where $B_j$ is the Blaschke product and $\alpha_j,\alpha_0,a_j$
fulfil
the relations \thetag{3}, \thetag{5} and \thetag{6}.\newline
Moreover, if $p_j\geq \frac12$ for some $j$, then we can
choose either
$B_j\equiv 1$ or
$B_j(\lambda)=\frac{\lambda-\alpha_j}{1-\bar\alpha_j\lambda}$ with
$|\alpha_j|<1$.\newline
Additionally, if $|\alpha_j|<1$ for all $j=1,\ldots,n$ then
either $B_j\equiv 1$ or
$B_j(\lambda)=\frac{\lambda-\alpha_j}{1-\bar\alpha_j\lambda}$ for
all
$j=1,\ldots,n$.
\endproclaim
Proposition 2 shows that $\kappa$-geodesics in $\Cal E(p)$ are,
in general, of a similar form as in the convex case. Observe that
we do not
exclude the case that the Blaschke products appearing in the
formulas
for geodesics may have more than one zero. Even in the case they
have only
one, we do not claim that this one equals $\alpha_j$. In
particular,
this is not clear for those $j$-th components with $|\alpha_j|=1$
and
$p_j<\frac12$. Nevertheless, our experience so far has led to the
conjecture
that there are no $\kappa$-geodesics at all with some
$|\alpha_j|=1$,
where $p_j<\frac12$.
\par As an easy consequence of Proposition 2 we get (cf.
\cite{Po}):

\proclaim{Corollary 3} Let $\phi:E\longrightarrow \Cal E(p)$ be a
$\kappa$-geodesic for $(\phi(0),\phi^{\prime}(0))$ with
$\phi^{\prime}(0)\neq 0$. Then
$\phi^*(\partial
E)\subset \partial \Cal E(p)$.
\endproclaim
As usual, $\phi^*$ denotes the boundary values of $\phi$.
\par In the paper \cite{BFKKMP} the authors delivered an effective
formula for the Kobayashi metric in the convex ellipsoids of type
$\Cal E(1,m)$ (i.e. for $m\geq \frac12$). In our note we shall
find the
formulas for the Kobayashi metric of the ellipsoids $\Cal E(1,m)$
with
$m<\frac12$. We obtain the formulas using Proposition 2,
proceeding similarly
as in \cite{JP} (we even use similar notation), where
the authors used Theorem 1 to get the formulas from \cite{BFKKMP}.

Let us make here one general remark: whenever, in the sequel, we
shall consider
$\Cal E(1,m)$, then we mean $m<\frac12$ (unless otherwise stated).

\par Remark that there are the following automorpisms
$$
\Cal E(1,m)\owns (z_1,z_2)\longrightarrow
\left(\frac{z_1-a}{1-\bar a
z_1},\frac{e^{i\theta}(1-|a|^2)^{\frac{1}{2m}}z_2}{(1-\bar a
z_1)^{\frac{1}{m}}}\right)\in \Cal E(1,m),
$$
where $a\in E$, $\theta\in \bold R$.
Therefore, in order to find the formulas for the Kobayashi metric,
it suffices to calculate the Kobayashi metric for $((0,b),(X,Y))$,
where
$b\geq 0$.
\par The cases: $b=0$, $X=0$ or $Y=0$ are as usual easily done.
\par In the case $b>0$, $XY\neq 0$ we shall need some more work
and much
more calculations. To make the calculations simpler and the
formulas for the
Kobayashi metric clearer we introduce in this case some additional
notation. Let us put
$$
v:=\left(\frac{b|X|}{m|Y|}\right)^2.
$$
Without loss of generality we may assume that $Y=1$.
For $v\leq \frac{1}{4m(1-m)}$ we define
$$
t:=\frac{2m^2v}{1+2m(m-1)v+\sqrt{1+4m(m-1)v}}.\tag{7}
$$
Observe that $\frac{1}{4m(1-m)}>1$ and
that $t$ increases when $v$ increases and, additionally, we have $t\leq
\frac{m}{1-m}<1$. Consequently, for any $v\leq \frac{1}{4m(1-m)}$
there is
exactly one solution (from the interval $(0,1)$) of the following
equation:
$$
x^{2m}-tx^{2m-2}-(1-t)b^{2m}=0.\tag{8}
$$
Then the formulas for the Kobayashi metric (we often write
$\kappa(v)$
instead of writing $\kappa_{\Cal
E(1,m)}((0,b);(X,Y))$) are given as follows:

\proclaim{Theorem 4} If $m<\frac12$ then the following formulas
hold:
$$
\align
\kappa_{\Cal E(1,m)}((0,0);(X,Y))&=h(X,Y),\\
\kappa_{\Cal E(1,m)}((0,b);(0,Y))&=\frac{|Y|}{1-b^2},\\
\kappa_{\Cal
E(1,m)}((0,b);(X,0))&=\frac{|X|}{(1-b^{2m})^{\frac12}},
\endalign
$$
where $h$ is the Minkowski function for $\Cal E(1,m)$.

And in the remaining cases ($b>0,X\neq 0,Y=1$)

\noindent (i) if $v\leq 1$, then
$$
\kappa(v)=
\frac{m}{b}\frac{x^{2m-1}}{(1-m)x^{2m}+mx^{2m-2}-b^{2m}}=:\kappa_1(v)\quad;
$$

\noindent (ii) if $v\geq \frac{1}{4m(1-m)}$, then
$$
\kappa(v)=\frac{m}{b}\frac{\sqrt{(1-b^{2m})v+b^{2m}}}{1-b^{2m}}=:\kappa_2(v)
\quad;
$$

\noindent (iii) if $1<v<\frac{1}{4m(1-m)}$, then
$$
\kappa(v)=\min\{\kappa_1(v),\kappa_2(v)\},
$$
where $x$ is the only solution in $(0,1)$ of equation \thetag{8}.

Moreover, in the formula (iii) the minimum is equal to
$\kappa_1(v)$ for $v\leq v_0$ and equal to $\kappa_2(v)$ for
$v>v_0$, where
$v_0:=\frac{t_0}{(t_0(1-m)+m)^2}$,
$t_0:=\frac{x_0^{2m}-b^{2m}}{x_0^{2m-2}-b^{2m}}$ and $x_0$ is the
only
solution in the interval $(0,1)$ of the equation
$$
\multline
x^{4m-2}(-1-2m+2m^2+b^{2m})+x^{2m}(1+(1-2m)b^{2m})+
x^{2m-2}(1+(2m-1)b^{2m})-\\
(1-m)^2x^{4m}-m^2x^{4m-4}-b^{2m}=0.
\endmultline
$$
\endproclaim

The detailed discussion of this situation also leads to the
following
properties, which
differ the non-convex case from the convex one (cf. \cite{BFKKMP},
\cite{M},
\cite{JPZ}).

\proclaim{Corollary 5} In the ellipsoid $\Cal E(1,m)$ the
following
properties hold:

\noindent (i) There are non-constant analytic discs of the form
\thetag{2}
with \thetag{3}--\thetag{6}, which are not $\kappa$-geodesics for
$(\phi(0),\phi^{\prime}(0))$.

\noindent (ii) There are non-constant analytic discs of the form
\thetag{2}
with \thetag{3}--\thetag{6}, which are not $k$-geodesics for
$(\phi(\lambda_0),\phi(\lambda_1))$ for some
$\lambda_0,\lambda_1\in E$
(see (i) and \cite{Ven}).

\noindent (iii) For any $b>0$ there is exactly one $v\in
(0,\infty)$
(equal to $v_0$ from Theorem 4) such that there are more than one
(up to a
M\"obius transformation, exactly
two) $\kappa$-geodesics with respect to the pair
$(b,v)$.

\noindent (iv) For any $b>0$ the function $\kappa_{\Cal
E(1,m)}((0,b);
(\cdot,1))$ is not differentiable at the point $X_0$ given by the
equation
$v_0=\left(\frac{b|X_0|}{m}\right)^2$.

\endproclaim
\par It seems that all the analytic discs of \thetag{2} with
\thetag{3}--\thetag{6} describe all the local $\kappa$-extremals.
\it In any case, our paper shows that there are a lot of local
$\kappa$-extremal and stationary
maps, which are not $\kappa$-extremals, even in the case of
strongly pseudoconvex Reinhardt domains, which arise from this
example
by smoothing the corners \rm (compare the examples of N. Sibony in
\cite{Pa}).

\subheading{2. Proof of Proposition 2} In the sequel we shall make
use of the
following characterization of $\kappa$-geodesics in a smooth
strongly
pseudoconvex domain $\Omega$ (writing smooth we always mean
$C^{\infty}$).
\par For a mapping $\phi\in \Cal O(E,\Omega)\cap C^1(\bar E,\bar
\Omega)$
with $\phi(\partial E)\subset \partial \Omega$, $\phi(0)\in
\Omega$ one
knows that
$$
\partial E\owns\lambda\longrightarrow \lambda \frac{\partial
\rho}{\partial z}
\left(\phi(\lambda)\right)\cdot\phi^\prime(\lambda)
$$
is a positive function (see \cite{Pa}),
where $\Omega=\{\rho<0\}$ and $\rho$ is a defining function of
$\Omega$
plurisubharmonic in
$\bold C^n$. Here $w\cdot z:=w_1z_1+\ldots+w_nz_n$ for
$w=(w_1,\ldots,w_n)$,
$z=(z_1,\ldots,z_n)\in \bold C^n$.
\par We define
$$
\gather
p(\lambda)^{-1}:=\lambda\frac{\partial \rho}{\partial
z}(\phi(\lambda))\cdot\phi^{\prime}(\lambda),\;\lambda\in \partial
E,\\
\tilde \phi
(\lambda):=p(\lambda)\lambda\frac{\partial\rho}{\partial z}
\left(\phi(\lambda)\right),\;\lambda\in\partial E.
\endgather
$$

\proclaim\nofrills{Definition 6}\;\rm{(see \cite{Pa},\cite{L}).} A
mapping
$\phi:\bar E\longrightarrow \bar \Omega$
is said to be \it stationary \rm if $\phi$ is a $C^1$ mapping of
$\bar E$ into $\bar \Omega$, holomorphic on $E$ such that
$\phi(\partial E)\subset\partial \Omega$, $\phi(0)\in \Omega$ and
$\tilde \phi$ extends to a continuous mapping on $\bar E$,
holomorphic on $E$.
\endproclaim
\noindent Observe that there are weaker notations of what
''stationary''
means (see \cite{JP}).
\par Then in \cite{Pa} the following result is shown:

\proclaim\nofrills{Theorem 7}\;{\rm (see \cite{Pa}).} Let
$\phi:\bar E\longrightarrow \bar \Omega$
be a $C^1$-mapping with $\phi(\partial E)\subset \partial \Omega$,
$\phi(0)\in E$, which is holomorphic on $E$. If $\phi$ is a
$\kappa_{\Omega}$-geodesic for
$(\phi(0),\phi^{\prime}(0))\in \Omega\times (\bold C^n\setminus
\{0\})$,
then $\phi$ is stationary.
\endproclaim

\par We point out that even more is shown in \cite{Pa}, namely,
that any
$\phi$ as in Theorem 7, which is a local extremal for $(z,X)$ ---
observe that
the term $\kappa$-geodesic means a global extremal w.r.t.
\thetag{1} ---
is automatically stationary.
\par It may be possible that one can find the proof of Proposition
2
(or even a better description of $\kappa$-geodesics than the one
given there)
by modifying slightly the solution of the extremal problem in
$C^1$
pseudoconvex domains given in \cite{Po}. Because of Corollary
5.(i) we are
not interested in this here.
\par We shall often make use of the decomposition theorem for a
mapping $f\in
H^{\infty}(E)$, $f\not \equiv 0$, which says that
$$
f(\lambda)=B(\lambda)A(\lambda)\text{ for $\lambda\in E$,}
$$
where $B$ is the Blaschke product of $f$ and $A$ is a nowhere
vanishing function
from $H^{\infty}(E)$.
\par Now we start the proof of Proposition 2 proving a sequence of
lemmas.

\proclaim{Lemma 8} Let $\phi:E\longrightarrow \Cal E(p)$ be a
$\kappa$-geodesic for $(\phi(\lambda_0),\phi^{\prime}(\lambda_0))$
with
$\phi^{\prime}(\lambda_0)\neq 0$ for $\lambda_0\in E$, where
$\phi_j\not \equiv 0$ for all $j$. Assume that
$$
\phi_j(\lambda)=B_j(\lambda)A_j(\lambda),\quad j=1,\ldots,n,
$$
where $B_j$ and $A_j$ are, as above, the factors from the
decomposition of $\phi_j$.
Denote by $\Cal Z_j$ the zeros of $B_j$ (counted with
multiplicity).
Denote by $\tilde\Cal Z_j$ a subset of $\Cal Z_j$. Let us
associate
with  $\tilde \Cal Z_j$ the Blaschke product $\tilde  B_j$. Put
$$
\tilde \phi_j:=\tilde B_jA_j\text{ and }
\tilde \phi:=(\tilde \phi_1,\ldots,\tilde \phi_n).
$$
If $\tilde \phi$ is non-constant, then $\tilde \phi$ is a
$\kappa$-geodesic for $(\tilde \phi (\lambda_0),\tilde
\phi^{\prime}(\lambda_0))$ and $\tilde\phi^{\prime}(\lambda_0)\neq
0$.
\endproclaim
\demo{Proof} At first we prove that $\tilde \phi(E)\subset \Cal
E(p)$. To see
this let us repeat inductively the following procedure. We divide
the
components $\phi_j$ by the Blaschke factor assigned to a zero from
$\Cal Z_j\setminus \tilde \Cal Z_j$ (or if we have already
exhausted
$\Cal Z_j\setminus \tilde \Cal Z_j$ we leave the component without
change).
We proceed so till we have exhausted all the sets $\Cal
Z_j\setminus \tilde
\Cal Z_j$
($j=1,\ldots,n$), so, in other words, till we get $\tilde \phi$.
Let us put
$\tilde h(z):=|z_1|^{2p_1}+\ldots+|z_n|^{2p_n}$.
Remark that in view of the maximum principle for subharmonic
functions applied to the composition of the mappings obtained from
$\phi$
after a finite number of steps of the above mentioned procedure
with $\tilde
h$ we get that this composition is not larger than $1$ on $E$.
Consequently the limit function of the composition is not larger
than $1$.
So we get $\tilde h\circ \tilde \phi\leq 1$ on $E$. The maximum
principle
implies that either $\tilde h\circ \tilde \phi\equiv 1$ or
$\tilde h\circ\tilde \phi<1$ on $E$. The first case gives, using
local peak functions, that $\tilde \phi$ is constant. This
completes
the proof of the inclusion above.
\par Suppose now that
$\tilde \phi$ is not a $\kappa$-geodesic for
$(\tilde \phi (\lambda_0),\tilde \phi^{\prime}(\lambda_0))$.
Then there exists a map $\tilde \psi:E\longrightarrow \Cal E(p)$
such that
$$
\tilde \psi (\lambda_0)=\tilde \phi(\lambda_0),\;
\tilde \psi ^{\prime}(\lambda_0)=\tilde
\phi^{\prime}(\lambda_0),\;
\tilde \psi(E)\subset \subset \Cal E(p).
$$
Let us define
$$
\psi_j(\lambda):=\frac{B_j(\lambda)}{\tilde B_j(\lambda)}
\tilde\psi_j(\lambda)\text{ for $\lambda\in E$ and }
\psi:=(\psi_1,\ldots,\psi_n).
$$
We have
$$
\psi(\lambda_0)=\phi(\lambda_0),\;\psi^{\prime}(\lambda_0)=
\phi^{\prime}(\lambda_0)
$$
(recall that $\phi_j(\lambda)=\frac{B_j(\lambda)}{\tilde
B_j(\lambda)}
\tilde\phi_j(\lambda)$) and $\psi(E)\subset\subset \Cal E(p)$.
But this contradicts the fact that $\phi$ is a $\kappa$-geodesic
for
$(\phi(\lambda_0),\phi^{\prime}(\lambda_0))$.
\par To prove that $\tilde \phi^{\prime}(\lambda_0)\neq 0$ we
proceed
by contradiction. Thus if $\tilde\phi^{\prime}(\lambda_0)=0$ then
$\eta:=(\eta_1,\ldots,\eta_n)$,
$\eta_j(\lambda):=\frac{B_j(\lambda)}{\tilde B_j(\lambda)}\tilde
\phi_j(\lambda_0)$, is also a $\kappa$-geodesic for
$(\phi(\lambda_0),\phi^{\prime}(\lambda_0))$. But $\eta(E)\subset
\subset
\Cal E(p)$ ---  a contradiction.
\qed\enddemo

\proclaim{Lemma 9} Let $\phi:E\longrightarrow \Cal E(p)$ be a
$\kappa$-geodesic for $(\phi(\lambda_0),\phi^{\prime}(\lambda_0))$
with
$\phi^{\prime}(\lambda_0)\neq 0$, where
$\phi_j\not\equiv 0$, $j=1,\ldots,n$. Then
$$
\phi_j(\lambda)=B_j(\lambda)\left(a_j\frac{1-\bar\alpha_j\lambda}
{1-\bar\alpha_0\lambda}\right)^{\frac{1}{p_j}},
$$
where $B_j$ is the Blaschke product and the coefficients
$\alpha_j,\alpha_0,a_j$ for $j=1,\ldots,n$ fulfil the relations
\thetag{3}, \thetag{5} and \thetag{6}.
\endproclaim
\demo{Proof} We know (from the decomposition theorem) that
$\phi_j=B_jA_j$, where $B_j$ is the Blaschke product and $A_j$ has
no zero
in $E$.
\par If $A_j$ is constant for $j=1,\ldots,n$, then we are done
with
$\alpha_0=\ldots=\alpha_n=0$ and $a_j=A_j^{p_j}$ because
$|B_j^*|=1$ a.e. on $\partial E$ implies
$\sum_{j=1}^{n}|A_j|^{2p_j}=1$
(otherwise, if the sum is smaller than $1$, then
$\phi(E)\subset\subset
\Cal E(p)$, hence $\phi$ is not a $\kappa$-geodesic).
So we may assume that some $A_j$ is not constant. For
$j=1,\ldots,n$ let us
put $\psi_j:=A_j$. Then in view of Lemma 8 the mapping
$$
\psi:=(\psi_1,\ldots,\psi_n)\text{ is a $\kappa$-geodesic for
$(\psi(\lambda_0),\psi^{\prime}(\lambda_0))$ with
$\psi^{\prime}(\lambda_0)\neq 0$.}\tag{9}
$$
Let us take $k\in \bold N$ \ such that $q_j:=p_jk\geq\frac12$
for $j=1,\ldots,n$.
Put $\tilde \psi_j:=\psi_j^{\frac1k}$, $q:=(q_1,\ldots,q_n)$,
$\tilde
\psi:=(\tilde \psi_1,\ldots,\tilde \psi_n)$.
Remark that $\tilde \psi (E)\subset \Cal E(q)$
and that $\Cal E(q)$ is a convex ellipsoid.
\par Now we will prove that $\tilde \psi$ is a $\kappa$-geodesic
for
$(\tilde \psi(\lambda_0),\tilde \psi^{\prime}(\lambda_0))$ in
$\Cal E(q)$
(with $\tilde\psi^{\prime}(\lambda_0)\neq 0$),
so; consequently, it is a complex geodesic.
\par To see this, remark that otherwise there would be
$\eta:E\longrightarrow\Cal E(q)$ such that
$\eta(\lambda_0)=\tilde \psi (\lambda_0)$,
$\eta^{\prime}(\lambda_0)=
\tilde \psi^{\prime}(\lambda_0)$ and $\eta(E)\subset\subset \Cal
E(q)$.
Hence $\eta^k:=(\eta_1^k,\ldots,\eta_n^k)$ maps $E$ into $\Cal
E(p)$
with $\eta^k(\lambda_0)=\psi(\lambda_0)$,
$(\eta^k)^{\prime}(\lambda_0)=\psi^{\prime}(\lambda_0)$ and
$\eta^k(E)\subset\subset \Cal E(p)$, which contradicts \thetag{9}.
\par Observe that $\tilde \psi_j$ is without zeros in $E$,
so in view of Theorem 1 we know that
$$
\tilde \psi_j(\lambda)=\left(a_j\frac{1-\bar\alpha_j\lambda}
{1-\bar\alpha_0\lambda}\right)^{\frac{1}{p_jk}},
$$
where $a_j,\alpha_j,\alpha_0$ fulfil the relations \thetag{3},
\thetag{5} and
\thetag{6}. This completes the proof.
\qed\enddemo

\proclaim{Lemma 10} Let $\phi:E\longrightarrow \Cal E(p)$ be a
$\kappa$-geodesic for $(\phi(\lambda_0),\phi^{\prime}(\lambda_0))$
with
$\phi^{\prime}(\lambda_0)\neq 0$, where
$\phi_j\not\equiv 0$ for $j=1,\ldots,n$. Assume
that $p_1,\ldots,p_k\geq\frac12$ for some $0\leq k\leq n$. Then
$B_j$, $j=1,\ldots,k$, from Lemma 9 can be chosen as
$\left(\frac{\lambda-\alpha_j}{1-\bar
\alpha_j\lambda}\right)^{r_j}$ with $r_j=0$ or $1$ (observe that the $a_j$'s
have to be modified by rotations); the
$\alpha_j$'s
are those from Lemma 9.
\endproclaim
\demo{Proof} From Lemma 9 we get
$$
\phi_j(\lambda)=B_j(\lambda)\left(a_j\frac{1-\bar\alpha_j\lambda}
{1-\bar\alpha_0\lambda}\right)^{\frac{1}{p_j}}
$$
for $j=1,\ldots,n$, where $B_j$ is the Blaschke product and the
coefficients
$\alpha_j,\alpha_0,a_j$ fulfil the relations \thetag{3},
\thetag{5} and \thetag{6}.
\par Put
$$
\psi_j:=\phi_j,j=1,\ldots,k,\quad\psi_j(\lambda):=
\left(a_j\frac{1-\bar\alpha_j\lambda}
{1-\bar\alpha_0\lambda}\right)^{\frac{1}{\ell
p_j}},j=k+1,\ldots,n,
$$
where $\ell\in\bold N$ is such that $\ell p_j\geq\frac12$,
$j=k+1,\ldots,n$.
If $\psi:=(\psi_1,\ldots,\psi_n)$ is constant then we are
done. Assume now that $\psi$ is
not constant. Then exactly as in the proof of Lemma 9 and in view
of Lemma 8
we get that $\psi$ is a complex geodesic for
$(\psi(\lambda_0),\psi^{\prime}(\lambda_0))$ with
$\psi^{\prime}(\lambda_0)\neq 0$ in
$\Cal E(p_1,\ldots,p_k,\ell p_{k+1},\ldots,\ell p_n)$ --- a convex
ellipsoid
--- which in view of Theorem 1 completes the proof.
\qed\enddemo
\demo{Proof of Proposition 2} If we combine Lemmas 9 and 10, then
we only have
to prove the last part of the proposition.
\par In view of Lemma 8 it suffices to discuss the case when the
Blaschke
products of all components $\phi_j$ of $\phi$ have at most a
finite number of
zeros.
\par Remark that in this case the mapping $\phi$ is holomorphic in
a
neighborhood of $\bar E$ and moreover (since the zeros of the
$\phi_j$'s are
lying in $E$ and their number is finite) it is
'far' from the points, where the boundary of the ellipsoid is not
strongly
pseudoconvex. Therefore it is a $\kappa$-geodesic in some smooth
strongly
pseudoconvex subdomain of $\Cal E(p)$, whose boundary coincides
with
the boundary of $\Cal E(p)$ everywhere except for a small
neighborhood of that part of $\partial \Cal E(p)$, where the
strong
pseudoconvexity breaks down. In particular, we may assume that the
defining
function of
the new domain coincides with the defining function of $\Cal E(p)$
($=|z_1|^{2p_1}+\ldots+|z_n|^{2p_n}-1$) near $\phi(\partial E)$.
Applying Theorem 7 we
get that $\phi$ is stationary. Therefore
$$
\frac{1}{\lambda}\tilde
\phi_j(\lambda)\phi_j(\lambda)=p(\lambda)p_j\left|\frac{1-\bar\alpha_j\lambda}{1
-\bar
\alpha_0\lambda}\right|^2|a_j|^2>0 \text{ on
$\partial E, j=1,\ldots,n$.}
$$
In view of Gentili's result (see \cite{G}) we obtain that
$$
\tilde\phi_j(\lambda)\phi_j(\lambda)=
r_j(\lambda-\gamma_j)(1-\bar\gamma_j\lambda)\text{ for $\lambda\in
\bar E,
r_j>0,\gamma_j\in\bar E$.}
$$
This implies that $\phi_j$ has at most one zero in $E$.
\par From the previous formulas we have:
$$
p(\lambda)p_j|a_j|^2\left|\frac{1-\bar\alpha_j\lambda}
{1-\bar\alpha_0\lambda}\right|^2=|\tilde\phi_j(\lambda)\phi_j(\lambda)|=
r_j|1-\bar\gamma_j\lambda|^2
$$
for $\lambda\in\partial E$.
So
$$
\frac{r_j|1-\bar
\gamma_j\lambda|^2}{p_j|a_j|^2|1-\bar\alpha_j\lambda|^2}=
\frac{p(\lambda)}{|1-\bar\alpha_0\lambda|^2},\text{ if }
\lambda\in\partial E.
$$
Consequently either
$$
\gamma_1=\dots=\gamma_n,\;\alpha_1=\dots=\alpha_n
$$
or
$$
\gamma_j=\alpha_j,\text{ for $j=1,\ldots,n$.}
$$
In the second case we are done. If we assume the first one, then
(from the
formulas \thetag{5} and \thetag{6})
$$
\gather
\alpha_0=\left(\sum_{j=1}^{n}|a_j|^2\right)\alpha_1,\\
1+|\alpha_0|^2=\left(\sum_{j=1}^{n}|a_j|^2\right)(1+|\alpha_1|^2).
\endgather
$$
It follows that
$$
\gather
\alpha_1(1+|\alpha_0|^2)=\alpha_0(1+|\alpha_1|^2),\text{ so }\\
\alpha_1-\alpha_0=\alpha_0\alpha_1(\bar\alpha_1-\bar\alpha_0).
\endgather
$$
If $\alpha_1\neq \alpha_0$, then
$|\alpha_0\alpha_1|=1$ --- a contradiction.
If $\alpha_1=\alpha_0$, then we have $\gamma_1=\ldots=\gamma_n\in
E$. We may
define $\alpha_j=\alpha_0:=\gamma_1$, which does not spoil earlier
relations and then we are done.
\qed
\enddemo

\subheading{3. Proof of Theorem 4} To get the formulas for the
Kobayashi
metric in the 'extremal' cases ($b=0$,
$X=0$ or $Y=0$) we need the following lemma.

\proclaim{Lemma 11} Put
$$
\phi:E\owns \lambda\longrightarrow
(\lambda^{r_1}z^0_1,\ldots,\lambda^{r_n}z^0_n)\in\Cal E(p),
$$
where $z^0=(z^0_1,\ldots,z^0_n)\in
\partial E(p)$, $r_j\in\{0,1\}$, $j=1,\ldots,n$
and $\#\{j:r_j=1\text{ and } z^0_j\neq 0\}\geq 1$. Then $\phi$
is a $\kappa$-geodesic for $(\phi(0),\phi^{\prime}(0))$.
\endproclaim
\demo{Proof} Without loss of generality we assume that
$$
\phi(\lambda)=(\lambda z^0_1,\ldots,\lambda z^0_k,
z^0_{k+1},\ldots,z^0_n),\tag{10}
$$
where $z^0_1\cdots z^0_k\neq 0$, $k+1\geq 2$, $z^0\in \partial
\Cal E(p)$.
\par Let us take a holomorphic mapping
$$
\gather
\psi:E\longrightarrow\Cal E(p)\text{ with
$\psi(0)=(0,\ldots,0,z^0_{k+1},\ldots,z^0_n)$
and}\tag{11}\\
\text{ $\psi^{\prime}(0)=t(z^0_1,\ldots,
z^0_k,0,\ldots,0)=t\phi^{\prime}(0),\;t>0$.}
\endgather
$$
Without loss of generality we may assume that $\psi$ is continuous
on
$\bar E$.
Let $\tilde h(z):=\sum_{j=1}^{n}|z_j|^{2p_j}$, $\tilde h\in PSH$.
In view of
\thetag{11} we may write
$$
\psi (\lambda)=(\lambda A_1(\lambda),\ldots,\lambda
A_k(\lambda),\psi_{k+1}(\lambda),\ldots,\psi_n(\lambda)).\tag{12}
$$
Put
$$
\tilde \psi:=(A_1,\ldots,A_k,\psi_{k+1},\ldots,\psi_n).
$$
Since $(\tilde h\circ\tilde \psi)(\lambda)\leq 1$ on $\partial E$,
we get that
$\tilde h\circ \tilde \psi\leq 1$ on $E$, so consequently
$\tilde h\circ \tilde \psi(0)\leq 1$ or
$$
\sum_{j=1}^{k}|A_j(0)|^{2p_j}+\sum_{j=k+1}^{n}|\psi_j(0)|^{2p_j}\leq
1.
$$
In view of \thetag{10}, \thetag{11} and \thetag{12} we see that
$$
\sum_{j=1}^{k}t^{2p_j}|z^0_j|^{2p_j}\leq
1-\sum_{j=k+1}^{n}|z^0_j|^{2p_j}=
\sum_{j=1}^{k}|z^0_j|^{2p_j}.
$$
Hence we obtain
$$
\sum_{j=1}^{k}|z^0_j|^{2p_j}(t^{2p_j}-1)\leq 0.
$$
But $p_j>0$, so $t\leq 1$, which completes the proof of the lemma.
\qed
\enddemo

\demo{Proof of Theorem 4} In view of Lemma 11 we get the formulas
in the
'extremal' cases. Therefore we may assume that $b>0$, $X\neq 0$
and $Y=1$.
\par Below we consider only the mappings of the following forms
$$
\gather
\phi(\lambda)=\left(\frac{a_1\lambda}
{1-\bar\alpha_0\lambda},\left(a_2\frac{1-\bar\alpha_2\lambda}
{1-\bar\alpha_0\lambda}\right)^{\frac1m}\right),\tag{13}\\
\phi(\lambda)=\left(\frac{a_1\lambda}
{1-\bar\alpha_0\lambda},\frac{\lambda-\alpha_2}{1-\bar\alpha_2\lambda}
\left(a_2\frac{1-\bar\alpha_2\lambda}
{1-\bar\alpha_0\lambda}\right)^{\frac1m}\right)\tag{14}
\endgather
$$
such that
$$
\phi(0)=(0,b),\tau\phi^{\prime}(0)=(X,Y)\tag{15}
$$
and $\alpha_j,\alpha_0,a_j$ fulfil \thetag{3}, \thetag{5} and
\thetag{6}.
\par In the sequel we shall find a mapping of form \thetag{13} or
\thetag{14} such that
$|\tau|$ from \thetag{15} is the smallest and, moreover, we shall
see that
$$
\gather
\text{the smallest value of $|\tau|$ in \thetag{15} is
never achieved by a mapping}\tag{16}\\
\text{of  type \thetag{13} with $|\alpha_2|=1$.}
\endgather
$$
\par We claim that this smallest
$|\tau|$ will be the value of the Kobayashi metric for \linebreak $((0,b),(X,Y))
$ (or
equivalently for suitable $v$). Why is it so? Theoretically the
smallest
$|\tau|$ may be provided by a mapping $\phi$ with a second
component such
that $|\alpha_2|=1$ whose Blaschke product has some zeros
(see Proposition 2
or Lemma 10). But if this were the
case then also the mapping with the 'deleted' Blaschke product of
the second
component would be a $\kappa$-geodesic for some other $v$ (see
Lemma 7). So
the new mapping would deliver the
minimum of the values of $|\tau|$ between the mappings of the
forms
\thetag{13} and \thetag{14} with \thetag{15} (for some other data
connected
with our new $v$), which however contradicts our earlier remark
(see \thetag{16}).
\par At this point we repeat, at least partially, the reasoning
from
\cite{JP}, where the authors calculate the Kobayashi metric for
convex ellipsoids $\Cal E(1,m)$.
\par Remark that there is a
mapping $\phi$ of form \thetag{13} iff $v\geq 1$ and in this case
$|\tau|$ is given by the formula (see \thetag{8.4.25} from
\cite{JP})
$$
|\tau|=\frac{m}{b}\frac{\sqrt{(1-b^{2m})v+b^{2m}}}{1-b^{2m}}.\tag{17}
$$
Moreover, there is a mapping $\phi$ as in \thetag{14} fulfilling
\thetag{15}
iff there is $0<x<1$ (and then $\alpha_1=x$) such that
$x^{2m-1}>b^{2m}$ (but
this holds always since $2m-1<0$ !) and
$$
v((m-1)x^{2m}-mx^{2m-2}+b^{2m})^2-
(x^{4m-2}-b^{2m}x^{2m}-b^{2m}x^{2m-2}+b^{4m})=0.\tag{18}
$$
And then our $|\tau|$ is given by the formula (one can easily
check that the
denominator in the formula below is positive).
$$
|\tau|=\frac{mx^{2m-1}}{b((1-m)x^{2m}+mx^{2m-2}-b^{2m})}\quad\text{
(see \cite{JP}
\thetag{8.4.32}).}
$$

\par Recall that the formula for $t$ (see \thetag{7}) has sense if
$v\leq
\frac{1}{4m(1-m)}$. We also remark that $v$ calculated from
\thetag{18}
fulfils the inequality $v\leq\frac{1}{4m(1-m)}$ (this implies that
there may
exist a mapping of the form \thetag{14} fulfilling \thetag{15}
only for
$v\leq\frac{1}{4m(1-m)}$);
moreover, we shall see that for any
$v\leq\frac{1}{4m(1-m)}$ there exists $x\in (0,1)$ fulfilling
\thetag{18} --- this will imply that in this case there is always
a mapping
of the form \thetag{14} with \thetag{15}.
\par The equality \thetag{18} is equivalent to
$$
\gather
x^{2m}-tx^{2m-2}-(1-t)b^{2m}=0\text{ or }\tag{19}\\
\text{ or }\\
(m-1)^2vx^{2m}-\frac{m^2v}{t}x^{2m-2}+\frac{1-v}{1-t}b^{2m}=0.\tag{20}
\endgather
$$
Remark that the equation \thetag{19} has exactly one solution
$x\in (0,1)$
for $0<v\leq\frac{1}{4m(1-m)}$ and that \thetag{20} has no
solution in the
interval $(0,1)$ for $0<v<1$ but has exactly one solution in
$(0,1)$ if
$1<v\leq\frac{1}{4m(1-m)}$ (for $v=1$ we have $x=1$ and this gives
a
mapping from \thetag{13} as above).
\par In order to choose a mapping of type \thetag{13} or
\thetag{14} with
\thetag{15} such that $|\tau|$ is minimal we have no problems for
$v<1$
(respectively $v>\frac{1}{4m(1-m)}$); these are the mappings of
type
\thetag{14} with $\alpha_1:=x$, where $x$ is a solution of
\thetag{19}
(respectively of type
\thetag{13)}. The problem is with $1\leq v\leq
\frac{1}{4m(1-m)}$ (or equivalently
$\left(\frac{m}{1-m}\right)^2\leq t\leq
\frac{m}{1-m}<1$). In this case we have to choose one of three (at
most)
mappings:
$$
\align
\tau_1(v)&=\frac{m}{b}\frac{x_1^{2m-1}}{(1-m)x_1^{2m}+mx_1^{2m-2}-b^{2m}};
\tag{21}\\
\tau_2(v)&=\frac{m}{b}\frac{\sqrt{(1-b^{2m})v+b^{2m}}}{1-b^{2m}};\tag{22}\\
\tau_3(v)&=\frac{m}{b}\frac{x_2^{2m-1}}{(1-m)x_2^{2m}+mx_1^{2m-2}-b^{2m}},
\tag{23}
\endalign
$$
where $x_1$ is the solution in $(0,1)$ of \thetag{19} and $x_2$ is
the
solution in $(0,1)$ of \thetag{20} (if $v=1$, then
$x_2=1$).
\par One can easily check that
$$
(m-1)^2vt^2-(1+2m(m-1)v)t+m^2v=0,
$$
so
$$
v=\frac{t}{(t(1-m)+m)^2}.\tag{24}
$$
Therefore we get from \thetag{20} (the equation w.r.t. $x_2$)
$$
(m-1)^2t(x_2^{2m}-b^{2m})=m^2(x_2^{2m-2}-b^{2m}).\tag{25}
$$
Substituting \thetag{19} into the formula \thetag{21} we may write
(we
consider now $\tau_j$'s as the functions of $t$)
$$
\tau_1(t)=\frac{m}{b}\frac{x_1^{2m-1}}{(x_1^{2m-2}-b^{2m})((1-m)t+m)}.
$$
Substituting \thetag{24} into \thetag{22} we get
$$
\tau_2(t)=
\frac{m}{b}\frac{\sqrt{(1-b^{2m})t+(t(1-m)+m)^2b^{2m}}}{(1-b^{2m})(t(1-m)+m)}.
$$
Substituting \thetag{25} into \thetag{23} we get
$$
\tau_3(t)=\frac{m}{b}\frac{x_2^{2m-1}}{(x_2^{2m-2}-b^{2m})}
\frac{(1-m)t}{m(m+(1-m)t)}.
$$
To obtain the formulas for the Kobayashi metric as in the theorem
it is enough
to prove that
$$
\tau_1(t)<\tau_3(t)\text{ for
$t\in\left(\left(\frac{m}{1-m}\right)^2,\frac{m}{1-m}\right)$}\tag{26}
$$
(for $t=\frac{m}{1-m}$ we have
$\tau_1(\frac{m}{1-m})=\tau_2(\frac{m}{1-m})$). In
particular, the inequality for $t=\left(\frac{m}{1-m}\right)^2$
(in other
words for $v=1$) will prove
that the mapping of type \thetag{13} with $|\alpha_2|=1$ is never
a
$\kappa$-geodesic (see \thetag{16}).

\par Therefore, to prove \thetag{26}, we calculate (using the
formulas for
$\tau_1$ and $\tau_3$ and obtaining $x_1^{\prime}(t)$ and
$x_2^{\prime}(t)$
from \thetag{19} and \thetag{25})
$$
\tau_1^{\prime}(t)-\tau_3^{\prime}(t)=
\frac{m}{2b}\frac{(m+t(m-1))}{(m+(1-m)t)^2}
\frac{(x_1x_2-1)(tx_2(1-m)-mx_1)}{(mx_1^2+t(1-m))(t(1-m)x_2^2+m)}.
$$
\par One sees that $m+t(m-1)>0$ and, after some calculations, that
$tx_2(1-m)-mx_1<0$ for $t\in
\left(\left(\frac{m}{1-m}\right)^2,\frac{m}{1-m}\right)$. This
implies that
$\tau_1^{\prime}(t)-\tau_3^{\prime}(t)$ is positive for these $t$.
Hence
\thetag{26} holds and this completes the proof of the formula for
the Kobayashi
metric.
\vskip2ex
\par To prove the last part of the theorem observe that it is
enough to
show (use the continuity of the Kobayashi metric of $\Cal E(1,m)$,
the
continuity of $\tau_1$ and $\tau_2$ and the fact that the
Kobayashi metric
equals $\tau_1(v)$ for $v\leq 1$ and $\tau_2(v)$ for
$v\geq\frac{1}{4m(1-m)}$) that there is only one $t_0\in
\left(\left(\frac{m}{1-m}\right)^2,\frac{m}{1-m}\right)$  such
that
$$
\gather
\tau_1(t_0)=\tau_2(t_0)\text{ or }\\
\frac{x_1^{2m-1}}{x_1^{2m-2}-b^{2m}}=
\frac{\sqrt{(1-b^{2m})t_0+(t_0(1-m)+m)^2b^{2m}}}{1-b^{2m}},\text{
where }\\
x_1^{2m}-t_0x_1^{2m-2}-(1-t_0)b^{2m}=0.
\endgather
$$
Substituting the second equality into the first one we end up with
the
following equality
$$
\multline
x_1^{4m-2}(-1-2m+2m^2+b^{2m})+x_1^{2m}(1+(1-2m)b^{2m})\\
+x_1^{2m-2}(1+(2m-1)b^{2m})-(1-m)^2x_1^{4m}-m^2x_1^{4m-4}-b^{2m}=0,
\endmultline
$$
which turns out to have exactly one solution $x_1\in(0,1)$. This
completes
the proof.
\qed\enddemo


\vskip4ex


\widestnumber\key {GGGGGHÙ}

\Refs


\ref \no BFKKMPÙ \by B.E.Blank, D.Fan, D.Klein, S.G.Krantz, D.Ma,
M.-Y.Pang
\paper The Kobayashi metric of a complex ellipsoid in $\bold C^2$
\jour Experimental
Math. \vol 1 \yr 1992 \pages 47-55
\endref
\ref \no GÙ \by G.Gentili \paper Regular complex geodesics in the
domain
$D_n=\{(z_1,\ldots,z_n)\in \bold C^n:|z_1|+\ldots+|z_n|<1\}$ \jour
Springer
Lecture Notes in Math. \vol 1275 \yr 1987 \pages 235-252 \endref
\ref \no JPÙ \by M.Jarnicki,P.Pflug \book Invariant Distances and
Metrics in
Complex Analysis \publ Walter de Gruyter \yr 1993 \endref
\ref \no JPZÙ \by M.Jarnicki,P.Pflug,R.Zeinstra \paper Geodesics
for convex
complex ellipsoids \jour Annali della Scuola Normale Superiore di Pisa
\vol XX
Fasc. 4 \yr 1993 \pages 535-543
\endref
\ref \no LÙ \by L.Lempert \paper La m\'etrique de Kobayashi et la
repr\'esentation des domaines sur la boule \jour Bull. Soc. Math.
France \vol
109 \yr 1981 \pages 427-479
\endref
\ref \no MÙ \by D.Ma \paper Smoothness of Kobayashi metric of
ellipsoids
\jour (preprint) \yr 1993
\endref
\ref \no PaÙ \by M.-Y.Pang \paper Smoothness of the Kobayashi
metric of
non-convex domains \jour International Journal of Mathematics \vol
4
\yr 1993 \pages 953-987
\endref
\ref \no PoÙ \by E.A.Poletskii \paper The Euler-Lagrange
equations
for extremal holomorphic mappings of the unit disk \jour Michigan
Math. J.
\vol 30 \yr 1983 \pages 317-333
\endref
\ref \no VenÙ \by S.Venturini \paper Intrinsic metrics
in complete circular domains \jour Math. Annalen \vol 288
\yr 1990 \pages 473-481
\endref
\ref \no VesÙ \by E.Vesentini \paper Complex geodesics
\jour Compositio Math. \vol 44 \yr 1981 \pages 375-394
\endref

\endRefs
\enddocument